\documentclass{amsart}

\usepackage{amsfonts,amssymb}
\usepackage{float}
\usepackage{hyperref}
\usepackage{enumerate}
\usepackage[capitalise]{cleveref}
\usepackage{amsmath, amsthm, mathabx}
\usepackage{verbatim}
\usepackage[official]{eurosym}
\usepackage{geometry}
\usepackage{nomencl}
\usepackage{bm}

\usepackage{fancyhdr}
\usepackage{amscd}
\usepackage[english]{babel}
\usepackage{xcolor}

\usepackage{todonotes}
\usepackage{lipsum}
\usepackage{amsthm}

\newtheorem{theorem}{Theorem}[section]
\newtheorem{proposition}[theorem]{Proposition}
\newtheorem{lemma}[theorem]{Lemma}
\newtheorem{corollary}[theorem]{Corollary}
\newtheorem{thmx}{Theorem}
\newtheorem{corx}[thmx]{Corollary}

\theoremstyle{definition}
\newtheorem{definition}[theorem]{Definition}

\DeclareMathOperator{\Sym}{Sym}

\DeclareMathOperator{\Cay}{Cay}
\DeclareMathOperator{\Sp}{Sp}
\DeclareMathOperator{\M}{M}

\title{Bounded conciseness in the space of marked groups}

\author{Federico Berlai}
\address[Federico Berlai]{Department of Mathematics, University of the Basque Country UPV/EHU, Barrio Sarriena s/n, 48940 Leioa, Spain}
\email[Federico Berlai]{federico.berlai@ehu.eus}

\date{\today}

\begin{document}

\begin{abstract}
We prove that bounded conciseness is a closed property in the space of marked groups. As a consequence, we reformulate a conjecture of Fern\'andez-Alcober and Shumyatsky \cite{FerSh} about conciseness in the class of residually finite groups.
\end{abstract}

\maketitle

\section{Introduction}
Conciseness draws its roots in the work of Philip Hall in the 50s. Given a word $w\in F(x_1,\ldots, x_n)$, one says that $w$ is concise in a group $G$ if, whenever the set of $w$-values $w\{G\}=\{w(g_1,\ldots, g_n)\mid g_i\in G\}$ is finite, then the subgroup $G_w$ generated by these finitely many values is also finite. A word is said to be concise in a class $\mathcal{C}$ of groups if it is concise in all groups of the class $\mathcal{C}$, and it is said to be concise if it is concise in the class of all groups. 
In the 50s, Hall asked if all words are concise. Since then many positive results have been obtained, but there are now expamples of non-concise words~\cite{Iv}. The similar question - but now restricted to the family of residually finite groups - is still open: that is, is every word concise in the class of residually finite groups?

Analogously, 
a word $w$ is \emph{boundedly concise} in a class of groups $\mathcal{C}$ if the following holds: there exists a function $\delta_{\mathcal{C},w}(m)=\delta(m)\colon \mathbb{N}\to \mathbb{N}$ such that whenever $G\in \mathcal{C}$ and $w\{G\}$ is finite, then $\lvert G_w\rvert \leqslant \delta\bigl( \lvert w\{G\}\rvert\bigr)$ (and in particular the group $G_w$ is finite).
By an ultraproduct argument \cite{FerMo}, it is known that if a word is concise in the class of all groups, then it is boundedly concise.


In this article,
after recalling known facts about the space of marked groups, we prove:

\begin{thmx}
Let $\mathcal{C}$ be a class of groups, $w\in F$ a word in a finitely generated free group, and suppose that the word $w$ is boundedly concise in $\mathcal{C}$. Then the word $w$ is boundedly concise in the class $\mathcal{L}$ of limits of groups in $\mathcal{C}$ in the space of marked groups.
\end{thmx}

LEF groups generalise residually finite groups, and are defined as follow. A group is said to be \emph{locally embeddable into finite groups} (abbreviated LEF) if for all finite subsets $F\subseteq G$ there exists a finite group $Q$ and a map $\varphi\colon G\to Q$ such that $\varphi\restriction_F$ is injective and $\varphi(f_1f_2)=\varphi(f_1)\varphi(f_2)$ for all $f_1,f_2\in F$.
If $G$ is residually finite, this map $\varphi$ is an actual homomorphism and $Q$ is a finite quotient of $G$. LEF groups are precisely the limits of finite groups in the space of marked groups.

\smallskip
From Theorem A we deduce:
\begin{corx}
Let $w\in F$ be a word. The following are equivalent:
\begin{enumerate}
    \item $w$ is boundedly concise in the class of residually finite groups;
    \item $w$ is boundedly concise in the class of LEF groups;
    \item $w$ is concise in the class of LEF groups.
    \end{enumerate}
\end{corx}
It is conjectured \cite{FerSh} that a word is concise in the class of residually finite groups if and only if it is boundedly concise in the same class. Thus, with Corollary B we can restate the conjecture as "a word is concise in the class of residually finite groups if and only if it is concise in the class of LEF groups".
Possibly, the latter reformulation can be of help in proving the conjecture.

\subsubsection*{Acknowledgements}
The author is supported by the Spanish Government, grant PID2020-117281GB-I00, partly by the European Regional Development Fund (ERDF), and the Basque Government, grant IT1483-22.
He is very thankful to Gustavo Fern\'andez, Matteo Pintonello, Andoni Zozaya for insightful conversations on this topic, and to an anonymous referee for their comments that improved a lot the quality of exposition of this work.

\section{Space of marked groups}\label{sec.spacemarked}
In this section, for the interested reader, we briefly discuss the space of marked groups, and we collect the proofs of some basic results that will be used later.

A \emph{marked group} is a pair $(G,S)$, where $G$ is a group and $S$ is a sequence of elements that generate $G$. 
If $\lvert S\rvert=n$ is finite, then $(G,S)$ is called an $n$-marked group.

If we have an $n$-marked group $(G,S)$, then there is a surjective homomorphism $\pi=\pi_{(G,S)}$ from the free group on $n$ free generators $x_1,\dots, x_n$ to the group $G$ defined on the generators as $\pi(x_i):=s_i$.  Two $n$-marked groups $(G,(s_1,\dots,s_n))$ and $(G',(s'_1,\dots,s'_n))$ are isomorphic
if the bijection $\varphi(s_i):=s_i'$ extends to an isomorphism of groups $\varphi\colon G\to G'$. In particular, two marked groups with given generating sets of different size are 
never isomorphic as marked groups, although they might be isomorphic as abstract groups.

Let $\mathcal{G}_n$ denote the set of $n$-marked groups, up to isomorphism of marked groups. Then $\mathcal{G}_n$ corresponds bijectively to the set of
normal subgroups of the free group $F_n$ on $n$ free generators $\{x_1,\dots, x_n\}$. It is a compact totally disconnected ultrametric space~\cite{ChGu}.

Let $(G,S),(G',S')\in\mathcal{G}_n$ and let $N$, $N'$ be the normal subgroups of $F=F_n$ such that $F/N\cong G$ and $F/N'\cong G'$. Let $B_{F}(r)$ denote the ball of radius $r$ in $F$ centered 
at the identity element with respect to the generating set $\{x_1^{\pm1},\dots, x_n^{\pm1}\}$. The function
\begin{equation}\label{definition_nu}
\nu(N,N'):=\sup\{r\in\mathbb N\mid N\cap B_{F}(r)=N'\cap B_{F}(r)\}\in \mathbb{N}\cup \{+\infty\}
\end{equation}
defines on $\mathcal{G}_n$ the ultrametric 
\begin{equation}\label{ultrametric}
d\bigl((G,S),(G',S')\bigr):=2^{-\nu(N,N')}.
\end{equation}
An $n$-marked group $(G,S)$ can be viewed as an $(n+1)$-marked group by adding the trivial element $e_G$ to $S$. This defines an isometric embedding of $\mathcal{G}_n$ into $\mathcal{G}_{n+1}$. Therefore, the ultrametrics given in Equation \eqref{ultrametric} extend to an ultrametric $d$ on $\mathcal{G}:=\bigcup_{n\in\mathbb N}\mathcal{G}_n$. This space is called the \emph{space of finitely generated marked groups}.



\medskip
Given a (marked) group $(G,S)$, one can define the Cayley graph $\Gamma:=\Cay(G,S)$, whose vertex set $V\Gamma$ is the group $G$ itself, and such that $(g,gs)$ is a directed edge whenever $g\in G$ and $s\in S\setminus \{e_G\}$. Thus, this graph is locally finite. In general $(gs,g)\notin E\Gamma$, unless $s^{-1}$ too belongs to $S$ (and not just $s$). 
There is also a natural labeling of the edges in $\Cay(G,S)$ by the set $S\setminus\{e_G\}$: any edge $(g,gs)$ is labeled by the letter $s$.

Given a marked group $(G,S)\in \mathcal{G}_n$, we can always assume $S$ to be a symmetric set by considering $(G,\tilde S)\in \mathcal{G}_{2n}$, where $\tilde S=S\cup S^{-1}$. When this happens, the Cayley graph $\Cay(G,\tilde S)$ can be viewed as an \emph{undirected} labeled graph, by replacing the edges $(g,gs)$, labeled by $s$, and $(gs,g)$, labeled by $s^{-1}$, with the undirected edge $\{g,gs\}$ labeled by~$s$. Nevertheless, we will not enforce this in what follows.

Let $(G,S)$ be a marked group. One can naturally define a norm $\lvert - \rvert_S\colon G\to \mathbb{N}$ in the following way:
\[\lvert g\rvert_S:=\min \{l\in \mathbb{N}\mid g=s_1\cdot\ldots\cdot s_l,\, s_i\in S\cup S^{-1}\,\forall\, i=1,\dots, l\}.\]
%
%
%
The following fact is well known. We provide a proof to make the paper more self-contained and accessible. 
The equivalence between the first and the third condition is discussed in \cite[Section 2.2]{ChGu}.

\begin{proposition}\label{folklore}
Let $(G,S)$ and $\{(G_r,S_r)\}_{r\in\mathbb N}$ be $n$-marked groups, 
and let $N$ and $\{N_r\}_{r\in\mathbb N}$ be normal subgroups of the free group $F$ on $n$ generators such that $F/N\cong G$ and ${F/N_r\cong G_r}$.
The following are equivalent:
\begin{enumerate}
\item the sequence $\{(G_r,S_r)\}_{r\in\mathbb N}$ converges to the point $(G,S)$ in the space $\bigl(\mathcal{G}_n,d\bigr)$;
\smallskip
\item for all $x\in N$ (respectively: for all $y\notin N$) there exists $\bar{r}$ such that $x\in N_r$ (respectively: $y\notin N_r$) for $r\geqslant\bar{r}$;
\smallskip
\item for all $R\in\mathbb N$ there exists $\bar{r}$ such that the Cayley graphs $\Cay(G,S)$ and $\Cay(G_r,S_r)$ have balls of radius $R$ isomorphic as labeled directed graphs, for all $r\geq\bar{r}$.
\end{enumerate}
\begin{proof}
Before starting with the proof we fix notation. Let $F=F(X)$ be the free group on $n$ free generators $X=\{x_1,\dots, x_n\}$, consider the surjective homomorfism $\pi\colon F\to G$ such that $\pi(x_i):=s_i$ for all $i=1,\dots, n$, where $S=(s_1,\dots, s_n)$, and the analogous surjective homomorfisms $\pi_r\colon F\to G_r$ such that $\pi_r(x_i):=s_{r,i}$, for all $i=1,\dots, n$, for all $r\in \mathbb{N}$, where $S_r=(s_{r,1},\dots, s_{r,n})$. By construction we have that $N=\ker (\pi)$ and that $N_r=\ker (\pi_r)$.
Let 
\[l_i:=\min\{l\in \mathbb{N}\mid s_i=wN \text{ and } \lvert w\rvert_X=l\}, \qquad \forall\, s_i\in S=(s_1,\dots, s_n),\]
that is, $l_i$ is the length in $F$ (with respect the word metric defined by the generating set~$X$) of a minimal representative of the coset $s_i$.
Set $L:=\max\{l_1,\dots, l_n\}$.

\smallskip
We now prove three implications. To begin, suppose that the sequence $\{(G_r,S_r)\}_{r\in\mathbb N}$ converges to the point $(G,S)\in \mathcal{G}_n$. In particular, for all $\varepsilon>0$ there exists $\delta\in \mathbb N$ for which \[d\bigl((G,S),(G_m,S_m)\bigr)=2^{ -\nu(N,N_m)} < \varepsilon,\qquad \forall\, m\geqslant \delta,\]
or, taking logarithms, such that $\nu(N,N_m)> - \log_2(\varepsilon)$. 
Let $x\in N$ be fixed (respectively, $y\notin N$ be fixed).
Therefore, there exists a minimal natural number $r\in \mathbb{N}$ such that $x\in B_{F}(r)$ (respectively, $y\notin B_{F}(r)$). Then, we can choose $\varepsilon_0$ so that $-\log_2(\varepsilon_0)> r$. Thus, with this choice of $\varepsilon_0$ (and its corresponding $\delta_0$), we have that $x\in N_m$ for all $m\geqslant \delta_0$ (respectively, $y\notin N_m$). We set $\bar r= \delta_0$.

\smallskip
Now, suppose that the second condition holds, and let $R\in \mathbb N$ be fixed.
Let us consider the finite set $B:=B_{F}(2 L R)\subseteq F$, the ball of radius $2LR$ (with respect to the word metric defined by $X$) centered at the identity of $F$. We apply condition $(2)$ to each of the finitely many elements in this ball and, after taking maximums, we obtain a natural number $\bar r$ such that for all elements $x\in B$, we have that $x\in N$ if and only if $x\in N_r$ for all $r\geqslant \bar r$. 

Let $ B_{G}(R)$ be the ball of radius $R$ centered at the identity in the Cayley graph $\Cay(G,S)$ and let $B_r$ be the ball of radius $R$ centered at the identity in the Cayley graph $\Cay(G_r,S_r)$. Let us define the maps $\varphi_r\colon B_{G}(R)\to B_r$, for all $r\geqslant \bar r$, in the following way:
if $g\in B_{G}(R)$ is a vertex, then $g=wN$ for some $w\in F$. Let $w$ be such an element of minimal length (for the choice of constants we made, we have that $w\in B_{F}(LR)$). Define $\varphi_r(g):=wN_r$.

First of all, let us check that $\varphi$ is a well-defined map. If $wN=w_0N $ for $w,w_0\in B$ of minimal length, then $w^{-1}w_0\in N$, and moreover $w^{-1}w_0\in B_{F}(2LR)$. Thus, it must be that $w^{-1}w_0\in N_r$ for all $r\geqslant \bar r$, and therefore $wN_r=w_0N_r$.

The map $\varphi_r$ is clearly surjective and injective.
Let us check that it preserves edges: if we consider two vertices $g, gs_i\in B_G(R)$ joined by a directed edge (labeled by $s_i$), for some $s_i\in S$, where $g=wN$ and $s_i=x_iN$, then $\varphi_r(g)=wN_r$ and
\[\varphi_r(gs_i)=\varphi_r (wx_i N)= (wx_i)N_r=wN_r \cdot x_iN_r=wN_r \cdot s_{r,i},\]
and therefore $(\varphi_r(g), \varphi_r(gs_i))$ is an edge in the Cayley graph $\Cay(G_r,S_r)$, labeled by $s_{r,i}$.
The argument is completely symmetric, because $x\in B\cap N$ if and only if $x\in B\cap N_r$. Thus, $\varphi_r$ is an isomorphism of labeled directed graphs, and condition $(3)$ is met.

\smallskip
To conclude, suppose that the third condition holds, and fix $\varepsilon > 0$. 
Given $r\in \mathbb{N}$, we have that $d\bigl((G,S),(G_r,S_r)\bigr)< \varepsilon$ if and only if $\nu(N,N_r)> -\log_2(\varepsilon)$. 
Notice, if the balls $B_G(R)$ and $B_r$ are isomorphic labeled graphs, then $\nu(N,N_r)\geqslant LR$.
Thus, if we take $R$ such that $LR> -\log_2(\varepsilon)$, then for the associated $\bar r$ we have that 
\[d\bigl((G,S),(G_r,S_r)\bigr)=2^{-\nu(N,N_r)}< \varepsilon, \qquad \forall\, r\geqslant \bar r,\]
and the proof is concluded.
\end{proof}
\end{proposition}

\section{Conciseness and bounded conciseness}
For the following definition we need to introduce some notation. Suppose that we are given a word $w$ in finitely many variables, that is an element of a finitely generated free group, say $w\in F_r$. This word defines an evaluation map
\[\begin{matrix} w\colon& G^r&\longrightarrow &G.\\
&(g_1,\dots,g_r)&\longmapsto & w(g_1,\dots,g_r)
\end{matrix}
\]
We denote by $w\{G\}$ the set of images of $w$ in the group $G$. Notice that this set may fail to be a subgroup: for instance, the set of commutators in a group $G$ (that is when we consider the word $w=[x,y]\in F_2$) in general is not a subgroup of $G$. We denote by $G_w$ the subgroup of $G$ generated by the set $w\{G\}$.

\begin{definition}
Let $w\in F_r$ be a word. We say that the word is \emph{concise} in a group $G$ if the following holds: whenever $w\{G\}$ is finite, then the group $G_w$ is also finite. Restricting this notion to a class of groups $\mathcal{C}$, one says that $w$ is concise in $\mathcal{C}$ if the above statement is met for all groups in the class $\mathcal{C}$.

A word $w$ in \emph{boundedly concise} in a class of groups $\mathcal{C}$ if the following holds: there exists a function $\delta_{\mathcal{C},w}(m)=\delta(m)\colon \mathbb{N}\to \mathbb{N}$ such that whenever $G\in \mathcal{C}$ and $w\{G\}$ is finite, then $\lvert G_w\rvert \leqslant \delta\bigl( \lvert w\{G\}\rvert\bigr)$ (and in particular the group $G_w$ is finite).
\end{definition}

Given a property of groups $\mathcal{P}$,
a group $G$ is said to be locally $\mathcal{P}$ if all finitely generated subgroups of $G$ satisfy property $\mathcal{P}$. For instance, a group is locally finite if all its finitely generated subgroups are finite, or it is locally residually finite if all its finitely generated subgroups are residually finite.

The following lemma is obvious.
\begin{lemma}\label{obvious_lemma}
Let $\mathcal{C}$ be a class of groups, $w\in F_r$ be a word, and suppose that $w$ is boundedly concise in $\mathcal{C}$. Then $w$ is boundedly concise in the class 
of groups that are locally $\mathcal{C}$,
maintaining the same function $\delta\colon \mathbb{N}\to \mathbb{N}$.
\end{lemma}

Thus, given a word $w\in F_r$ and a class of groups $\mathcal{C}$ in which $w$ is (boundedly) concise, in view of Lemma \ref{obvious_lemma} we can always assume that the class $\mathcal{C}$ is closed with respect to taking groups that are locally in $\mathcal{C}$.

It is not clear if the class of groups on which a word is concise is closed under taking subgroups. Nevertheless, the next weaker lemma will be useful when considering non-finitely generated groups.

\begin{lemma}\label{lemma_to_fg_subgroups}
Let $G$ be a group and $w\in F_r$ be a word such that $w\{G\}$ is a finite set. There exists a finitely generated subgroup $H\leqslant G$ such that $w$ is concise in $G$ if and only if $w$ is concise in $H$.
\begin{proof}
If $w\{G\}$ is a finite set, then there exists a finite subset $F$ of $G$ such that $w\{F\}=w\{G\}$. If we consider the finitely generated subgroup $H:=\langle F\rangle$ of $G$, we see that $w\{H\}=w\{G\}$. Thus, the claim follows.
\end{proof}
\end{lemma}

\begin{theorem}\label{thm_main}
Let $\mathcal{C}$ be a class of groups closed under subgroups, $w\in F$ a word in a finitely generated free group, and suppose that the word $w$ is boundedly concise in $\mathcal{C}$. Then the word $w$ is boundedly concise in the class $\mathcal{L}$ of limits of groups in $\mathcal{C}$ in the space of marked groups.

\begin{proof}

We start by assuming that $G$ is finitely generated, to then consider the general case.
Let $S$ be a finite generating set for $G$, and suppose that $(G,S)$ a limit of $\mathcal{C}$-groups, say limit of the sequence $\{(G_r,S_r)\}_{r\in\mathbb N}\subseteq \mathcal G$. If $w\{G\}$ is an infinite set, there is nothing to prove. Thus, suppose that $w\{G\}\subseteq G$ is finite, so that there exists $R\in \mathbb{N}$ such that $w\{G\}\subseteq B_{(G,S)}(R)$.

In view of Proposition \ref{folklore}, there exists $\bar r\in \mathbb{N}$ such that $\Cay(G,S)$ and $\Cay(G_r,S_r)$, for all $r\geqslant \bar r$, have balls of radius $R$ isomorphic as labeled graphs. In particular, $w\{G\}$ and $w\{G_r\}$ are identified under these isomorphisms, and thus $\lvert w\{G\}\rvert=\lvert w\{G_r\}\rvert$ for all $r\geqslant \bar r$.

As $G_r\in \mathcal{C}$ and $w\{G_r\}$ is a finite set, we have that $G_w$ is a finite group, because the word $w$ is (boundedly) concise in $\mathcal{C}$ and $G_r\in \mathcal{C}$. Moreover, as the word $w$ is boundedly concise in this class and the cardinality $\lvert w\{G_r\}\rvert=\lvert w\{G\}\rvert$ remains constant for all $r\geqslant \bar r$, we obtain a uniform bound $\delta=\delta\bigl(\lvert w\{G\} \rvert \bigr)\in \mathbb{N}$ such that the cardinality of the finite groups $ (G_r)_w$ is bounded by $\delta$, for all $r\geqslant \bar r$.

Replacing $R$ with a bigger constant, if needed, we can suppose without loss of generality that the balls of radius $R$ in $G_r$ contain the finite group $(G_r)_w$ (of cardinality at most $\delta$, uniformly), for all sufficiently large indices $r$, say for all $r\geqslant \tilde r\geqslant \bar r$. Thus, $G_w$ must be a finite subgroup of $G$, of cardinality bounded by $\delta$, because $G$ and $G_r$ have isomorphic balls of radius $R$ for all sufficiently large~$r$.

If $G$ is not finitely generated and $w\{G\}$ is a finite set (say, realised by a finite set $F\subseteq G$ of elements, that is $w\{G\}=w\{F\}$), then Lemma \ref{lemma_to_fg_subgroups} and \cite[Proposition~2.20]{ChGu} imply that the finitely generated subgroup $\langle F\rangle$ is a limit of $(H_r,T_r)$, for sufficiently large indices, where $H_r\leqslant G_r$ are given by Lemma \ref{lemma_to_fg_subgroups}. As $\mathcal{C}$ is closed under taking subgroups, we deduce that $H_r\in \mathcal{C}$ for all $r$. In particular, $(H_{r})_w$ is a finite subgroup of $H_r$, with uniform bound on the cardinality not depending on $r$, but just on $\lvert w\{F\}\rvert $, and we can apply what we proved in the previous part of the proof.
%
\end{proof}
\end{theorem}



This point of view is general, and does not depend on a specific class. Applying the previous theorem to the class of (residually) finite groups, we extend all known results on bounded conciseness for those groups to LEF groups (that is, limit of finite groups in the space of marked groups).

\begin{definition}\label{definitionLEF}
A group $G$ is \emph{locally embeddable into finite groups} (for short, a LEF group) if for all $F\subseteq G$ finite subset there exists a finite group $Q$ and a map $\varphi\colon G\to Q$ such that
\begin{enumerate}
    \item $\varphi\restriction_F$ is an injective map;
    \item $\varphi(gh)=\varphi(g)\varphi(h)$ for all $g,h\in F$.
\end{enumerate}
\end{definition}
The second condition of a LEF group guarantees that the map $\varphi$ behaves like a homomorphism on the finite set $F$, whereas the first condition excludes non-meaningful ways to do that (for instance, by mapping all elements of $F$ to the trivial element $e_Q$).

When the group $G$ is residually finite, then we can construct \emph{homomorphisms} onto finite quotients of $G$ with these properties. It turns out that, for finitely \emph{presented} LEF groups, that is always the case, that is finitely presented LEF groups are residually finite. This is not the case for finitely generated groups. For instance, the groups of permutations $S:=\Sym_f(\mathbb{Z})$ with finite support of an infinite set is an infinite non-residually finite group because it contains the simple infinite group of even permutations with finite support. Let us consider the translation $t\colon \mathbb{Z}\to \mathbb{Z}$ defined by $t(n):=n+1$ for all integer numbers $n$. Then, the group $S\rtimes \langle t\rangle$ is a finitely generated LEF group that is not residually finite (it is generated by $t$ and the transposition $(0\ 1)\in S$).

Notice that a limit of residually finite groups is LEF. Indeed, if $G$ is a limit of residually finite groups and $F$ is a finite subset of $G$, then there exists a radius $R$ such that the finite subset $F\cdot F=\{fh\mid f,h\in F\}$ is contained in the ball in $G$ of radius $R$. This finite ball is eventually isomorphic as a labeled directed graph to the ball of radius $R$ in a residually finite group $G_r$. Call $\psi$ the restriction of this graph isomorphism to the vertex set and extend arbitrarily $\psi$ to a map $G\to G_r$. As $G_r$ is residually finite, we can find a homomorphism $\pi\colon G_r\to Q$ onto a finite group $Q$ that is injective on the finite set $\psi(F\cdot F)$. Then $\pi\circ \psi\colon G\to Q$ is a map that meets the conditions of Definition \ref{definitionLEF} for the set $F$, proving therefore that $G$ is LEF.

As an application of Theorem \ref{thm_main} we obtain:

\begin{corollary}
Let $w\in F$ be a word. The following are equivalent:
\begin{enumerate}
    \item $w$ is boundedly concise in the class of residually finite groups;
    \item $w$ is boundedly concise in the class of LEF groups;
    \item $w$ is concise in the class of LEF groups.
\end{enumerate}
\begin{proof}
The implication $(1)\Rightarrow (2)$ is Theorem \ref{thm_main}, whereas $(2)\Rightarrow (1)$ and $(2)\Rightarrow (3)$ are clear from the definitions. The remaining implication follows from the fact that an untraproduct of LEF groups is itself a LEF group (this is mentioned/proved in \cite[Proposition 2.2]{OuldP}, but compare also the appendix of~\cite{FerMo}).
\end{proof}
\end{corollary}

\section{Some examples}
In this concluding section we collect some examples of groups that are far from being residually finite (thus, neither linear) such that all words are concise in these groups.

We record the following lemma, whose proof can be adapted from \cite[Proposition~1.1.3]{Seg}:

\begin{lemma}\label{virtuallyC}
Let $w\in F_r$ be a word and suppose that $w$ is concise in a class of groups $\mathcal{C}$. Then, it is concise in the class of finite-by-$\mathcal{C}$ groups.
\end{lemma}

Lemma \ref{virtuallyC} is of interest because there exist several finitely presented finite-by-linear (central) extensions that are very far away from being residually finite. One such family is given by quotients of Deligne's group \cite{CoGuPi,De78,Morris}. Deligne's group $D$ is a central extension $\{e\}\to \mathbb{Z}\to D\to \Sp_{2n}(\mathbb{Z})\to \{e\}$, where $\Sp_{2n}(\mathbb{Z})$ is the group of symplectic matrices with integer coefficients, that is
\[\Sp_{2n}(\mathbb{Z}):= \{P\in \M_{2n}(\mathbb{Z})\mid P^t\left(\begin{smallmatrix}0 & I_n\\
-I_n & 0 \end{smallmatrix}\right)P=\left(\begin{smallmatrix}0 & I_n\\
-I_n & 0 \end{smallmatrix}\right)\}.\]
Thus, in particular, $\Sp_{2n}(\mathbb{Z})$ is a linear group. The group $D$ is finitely presented because both $\mathbb{Z}$ and $\Sp_{2n}(\mathbb{Z})$ are finitely presented. What is relevant to us is that Deligne \cite{De78} proved that $D$ is not residually finite, because the intersection of all subgroups of finite index is a cyclic subgroup of order two. Moreover, it satisfies Kazhdan's property (T) \cite[Example~1.7.13 (iii)]{BeHaVa}, and in particular any amenable quotient of $D$ must be finite. Thus, $D$ is not residually amenable, nor a limit of amenable groups in the space of marked groups (being finitely presented, these two notions coincide).

\smallskip
This, of course, is a central extension with \emph{infinite} normal subgroup and linear quotient, and thus the hipoteses of Lemma \ref{virtuallyC} are not met. Nevertheless, we can consider the quotients $D/m\mathbb{Z}$ given by the central subgroups $m\mathbb{Z}$: these are central extensions of the form
\[\{e\}\to \mathbb{Z}/m\mathbb{Z}\to D/m\mathbb{Z}\to \Sp_{2n}(\mathbb{Z})\to \{e\},\]
thus finite-by-linear central extensions, and therefore any word $w$ is concise in these $D/m\mathbb{Z}$. Notice that these groups too are not residually finite, and they are isolated points in the space of marked groups \cite[Section 5.8]{CoGuPi}. 
All words are concise in the class of linear groups~\cite{Mer}. Therefore, by Lemma \ref{virtuallyC}, all words are concise in these (non-linear) examples.

\end{document}